\documentclass[12pt]{amsart}

\usepackage{amssymb,amsfonts,amsmath,graphicx}
\usepackage{tikz-cd}
\usepackage{enumerate}
\usepackage{mathrsfs, mathtools}
\usepackage{bbm}
\usepackage[margin=1in]{geometry}
\definecolor{brown(web)}{rgb}{0.65, 0.16, 0.16}
\definecolor{cred}{rgb}{0.7, 0.11, 0.11}
\usepackage{graphicx,subcaption}
\usepackage[colorlinks=true, allcolors=blue]{hyperref}
\usepackage[noabbrev,nameinlink,capitalise]{cleveref}
\usepackage{thmtools} 
\usepackage{cite}
\usepackage{quiver}

\newtheorem{thm}{Theorem}[section]
\newtheorem{cor}[thm]{Corollary}
\newtheorem{prop}[thm]{Proposition}
\newtheorem{lem}[thm]{Lemma}

\theoremstyle{definition}
\newtheorem{defn}[thm]{Definition}

\newtheorem{exmp}[thm]{Example}

\theoremstyle{remark}
\newtheorem{rem}[thm]{Remark}

\newcommand{\define}[1]{\textbf{#1}}

\usepackage{appendix}

\newcommand{\bA}{\mathbb{A}}
\newcommand{\bR}{\mathbb{R}}
\newcommand{\R}{\mathbb{R}}
\newcommand{\bS}{\mathbb{S}}
\newcommand{\bD}{\mathbb{D}}
\newcommand{\bX}{\mathbb{X}}
\newcommand{\pht}{\text{PHT}}
\newcommand{\PHT}{\text{PHT}}

\usepackage{soul}
\newcommand{\Vect}{\text{Vect}}
\newcommand{\bsd}{\mathbb{S}^{d-1}}
\newcommand{\bdd}{\mathbb{D}^{d}}
\newcommand{\CS}{\mathcal{CS}}
\newcommand{\Dgm}{\mathbf{Dgm}}

\newcommand{\bT}{\mathbb{T}}

\newcommand{\AX}{\mathbb{AX}}

\newcommand{\ECT}{\text{ECT}}
\newcommand{\BCT}{\text{BCT}}

\newcommand{\cC}{\mathcal{C}}
\newcommand{\cE}{\mathcal{E}}
\newcommand{\cD}{\mathcal{D}}
\newcommand{\bI}{\mathbb{I}}


\numberwithin{equation}{section}

\title[Kan Approximations of the PHT]{Kan Approximations\\ of the Persistent Homology Transform}
\author[S.~Arya]{Shreya Arya}
\address{Department of Mathematics, University of Pennsylvania, Philadelphia, PA, USA}
\email{smarya@sas.penn.edu}
\author[J.~Curry]{Justin Curry}
\address{Department of Mathematics and Statistics, University at Albany (SUNY), Albany, NY}
\email{jmcurry@albany.edu}

\dedicatory{Dedicated to the memory of Sayan Mukherjee (1971--2025)}
\date{\today}

\begin{document}

\begin{abstract}
    The persistent homology transform (PHT) of a subset $M \subset \mathbb{R}^d$ is a map $\text{PHT}(M):\mathbb{S}^{d-1} \to \mathbf{Dgm}$ from the unit sphere to the space of persistence diagrams. This map assigns to each direction $v\in \mathbb{S}^{d-1}$ the persistent homology of the filtration of $M$ in direction $v$. In practice, one can only sample the map $\text{PHT}(M)$ at a finite set of directions $A \subset \mathbb{S}^{d-1}$. This suggests two natural questions: (1) Can we interpolate the PHT from this finite sample of directions to the entire sphere? If so, (2) can we prove that the resulting interpolation is close to the true PHT? In this paper we show that if we can sample the PHT at the module level, where we have information about how homology from each direction interacts, a ready-made interpolation theory due to Bubenik, de Silva, and Nanda using Kan extensions can answer both of these questions in the affirmative. A close inspection of those techniques shows that we can infer the PHT from a finite sample of heights from each direction as well. Our paper presents the first known results for approximating the PHT from finite directional and scalar data. 
\end{abstract}

\maketitle
\section{Introduction}

The Euler Characteristic and Persistent Homology Transforms \cite{turner}, abbreviated ECT and PHT respectively, are directional topological transforms reminiscent of the Radon Transform in integral geometry.
First introduced in 1917 \cite{radon1917bestimmung,radon1986determination}, the Radon transform takes a function $f$ on $\bR^d$ and returns a function $\hat{f}$ on the Grassmanian of $(d-1)$-planes, given by integrating $f$ over each (affine) hyperplane $W$.
This style of transform was then generalized in the 1980s \cite{schapira1988cycles,viro1989some} and 1990s \cite{schapira1991operations,schapira1995tomography} in a more topological direction, by replacing integration of $f|_W$ with Euler characteristic $\chi$ of the intersection $M\cap W$, where $M\subseteq \bR^d$ is some suitably nice shape.
Although the original motivation for this Euler calculus \cite{curry2012euler} is most reliably traced to the study of characteristic classes \cite{macpherson1974chern} and index theory \cite{kashiwara1985index}, Schapira's early work \cite{schapira1991operations} also describes direct motivation by problems in computational geometry and robotics \cite{guibas1983kinetic}.
With the advent of persistent homology in the early 2000s \cite{edelsbrunner2002topological},
this trend towards applications has exploded with \cite{robinson2014topological,hofer2017deep,crawford2020predicting,amezquita2022measuring,roelldifferentiable,marsh2024detecting} as representative examples of applied research using topological transforms; see \cite{munch2023invitation} for a recent and approachable introduction to some of these ideas.

One major drawback of topological transforms is that interpolation and approximation is inherently difficult.
At first glance, the discrete nature of topological invariants seems to provide an insurmountable obstacle to interpolation.
For example, if a shape $M$ has $\chi(M\cap W)=3$ and $\chi(M\cap W')=4$ for two nearby hyperplanes $W$ and $W'$ it makes no sense to linearly interpolate or average these values for a third hyperplane $W''$ because the Euler characteristic $\chi$ is inherently integer-valued. 
Modern applied topology has learned to sidestep these kinds of apparent instabilities by working with 1-parameter families of topological invariants---suitably called \emph{persistent} topological invariants---where continuity between these signatures can be proved.
For example, if one filters a shape $M$ in direction $v$ and studies either the Euler curve or the Betti curves (one for each $n\in\{0,\ldots,d-1\}$)
\[
    \ECT(M)(v, t):=\chi(h_v(-\infty,t]) \qquad \text{or} \qquad \BCT_n(M)(v,t):=\beta_n(h_v(-\infty,t]),
\]
where $h_v(x):=x\cdot v$ is the height function in direction $v$, then \cite{curry_howmany} shows that all of these curves are continuous in the argument $v\in \bsd$, for target $L^p(\R)$ and $p\in[1,\infty)$.
However, once again, averaging or interpolating these curves between nearby directions does not make sense due to their piece-wise constant integer nature.
Surprisingly, another shift in perspective, where one presents the Euler and Betti curves as sums of carefully defined indicator functions---indicator functions supported on the bars in the persistent homology barcode---allows a continuous interpolation theory, explained below.
This also serves as a pedagogical point for demonstrating why the PHT's extra structure provides a distinct advantage over the ECT, even though both are equally expressive\footnote{``Expressive'' should be taken to mean expressive in it's ability to discriminate shapes. As \cite{ghrist2018persistent,curry_howmany} both show, the ECT and PHT are both injective, i.e., if $\ECT(M)=\ECT(N)$, then $M=N$.}.

As reviewed in \cref{sec:pers-modules-review} and \cref{sec:PHT-review}, the PHT of a shape $M$ in direction $v$ and homological degree $n$ is the persistence diagram obtained by filtering $M$ in direction $v$, i.e., it is the unique multi-set of intervals $\mathcal{D}_n(v)=\{[b_i,d_i)\}$ with the property that 
the rank of the map on homology $H_n(h_v(-\infty,t]) \to H_n(h_v(-\infty,s])$ can be computed by counting the intervals in $\mathcal{D}_n(v)$ containing $[s,t]$.
Specializing this property to $s=t$, and recalling that the Betti number $\beta_n$ is the dimension of the homology group $H_n$, shows that the Betti curve can be expressed as a sum of indicator functions supported on intervals in the diagram $\mathcal{D}_n(v)$:
\[
    \BCT_n(v,\bullet) = \sum_{[b_i,d_i)\in \mathcal{D}_n(v)} \mathbbm{1}_{[b_i,d_i)} (\bullet)
\]
This perspective suggests a path forward to interpolating topological information between directions. 
If one is given a partial matching between diagrams $\mathcal{D}_n(v)$ and $\mathcal{D}_n(v')$ for nearby directions $v$ and $v'$, then the straight line interpolation $(1-\lambda)(b_i,d_i) + \lambda (b'_i,d'_i)$ between the endpoints of each pair of matched intervals provides a path between these persistence diagrams, and hence a candidate extension of the PHT along a path between $v$ and $v'$.
Unfortunately, the data of a matching between persistence diagrams is often missing and even when one can be found, it might not be unique.
This lack of canonical matchings is at the heart of several subtle technical difficulties that one encounters when trying to perform regression in the space of persistence diagrams $\Dgm$:
\cite{mileyko2011probability} shows that sets of persistence diagrams do not have unique Fr\'echet means and \cite{munch2015probabilistic} shows that even when these means are unique they might not vary continuously.

In this paper we circumvent these problems and answer two fundamental questions about the PHT, both in the affirmative:
\begin{enumerate}
    \item (Interpolation) Given a set of directions $A \subset \mathbb{S}^{d-1}$ and $\PHT(M)$ for these directions, does there exist a continuous extension $\overline{\pht}(M): \mathbb{S}^{d-1} \to \Dgm$ of this information?
    \item (Approximation) Given an extension $\overline{\pht}(M)$, can we prove reasonable error bounds between it and the true $\PHT(M)$?
\end{enumerate}
We achieve these results by utilizing pre-existing work \cite{de2013geometry,bubenik} by Bubenik, de Silva, and Nanda, which shows how category theory can provide canonical methods---Kan extensions---for interpolating between persistence diagrams, \emph{as long as one works on the module level}. 
These methods are especially nice because the provided extensions are Lipschitz, thus ensuring that small changes in direction correspond to small changes in persistent homology.
In answering the second question, we prove that the distance between the true PHT of shape $M$ and the extension $\overline{\pht}(M)$ is bounded by twice the Hausdorff distance between the sphere and the approximating direction set $A$. 
This has a nice online learning interpretation: as persistent homology from more directions is learned, we can update our inferred $\overline{\pht}(M)$ to more closely approximate the truth.

We note that \cite{de2013geometry,bubenik} also provides some flexibility when choosing a metric on the sphere of directions.
There are two obvious choices here---the Euclidean metric $d_2$ and the round (or geodesic) metric $d_g$---and many of our results are stated for both choices.
The choice of the round metric achieves optimal error bounds between the Kan-extended $\overline{\pht}(M)$ and the true PHT, but this optimality comes at a conceptual cost.
Because the round metric is best understood in terms of an angular coordinate system on the sphere, our persistence modules are naturally indexed by the totally ordered set $\Theta=[-\pi,0]$, instead of $\R$. 
This requires invoking the language of generalized persistence modules \cite{bubenik_generalisedPM}, which introduces some subtleties when defining the interleaving distance.
On the other hand, working with the Euclidean metric trades optimal error bounds for conceptual simplicity and better matches existing literature on the PHT.
For this reason, we choose to present both using the notation $\PHT^X$ to indicate that the indexing set $X$ can change, depending on one's preference.
For experts tempted to skip \cref{sec:Background}, we caution them to pay careful attention starting at \cref{defn:map version2}, as our $\Theta$-coordinate system for the PHT is original to this work. 

In \cref{sec: lipschitz extension}, we demonstrate how the theory of Kan extensions provides a constructive answer to our first question (on interpolation) and compute the left, right, and center Kan extensions for a few simple examples.
In \cref{sec: error bounds}, we answer our second question (on the approximation error) by deriving precise error bounds.
We first prove a point-wise, direction-by-direction interleaving bound, between the true PHT and its extension and provide an example showing that our bound is tight.
This then leads naturally to a global error bound between the PHT and it's approximation in the $L^{\infty}$ (sup) norm.
Finally, we generalize our results to the most practical setting by providing an error bound for extensions from \emph{fully discrete data} where both directions and filtration values are sampled.



\section{Background}
\label{sec:Background}


In this section, we establish the mathematical framework necessary for our main results.

\subsection{O-minimality}
We work within the framework of o-minimal structures due to its tame topological properties. This framework encompasses a broad class of geometric objects while excluding pathological constructions such as the Cantor set.

\begin{defn}\label{defn:o-minimal}
An \define{o-minimal structure} $\mathcal{O}= \{\mathcal{O}_d\}$, is a specification of a boolean algebra of subsets $\mathcal{O}_d$ of $\bR^d$ for each natural number $d\geq 0$, satisfying certain axioms.
These are: we assume that $\mathcal{O}_1$ contains only finite unions of points and intervals.
We further require that $\mathcal{O}$ be closed under certain product and projection operators, i.e.~if $A \in \mathcal{O}_d$, then $A\times \bR$ and $\bR\times A$ are in $\mathcal{O}_{d+1}$; 
and if $A \in \mathcal{O}_{d+1}$, then $\pi(A)\in \mathcal{O}_d$ where $\pi: \bR^{d+1}\to \bR^d$ is axis-aligned projection. 
It is a fact that $\mathcal{O}$ contains all semi-algebraic sets, but may contain certain regular expansions of these sets.
Elements of $\mathcal{O}$ are called \define{definable sets} and \emph{compact} definable sets are called \define{constructible sets}.
The subcollection of constructible subsets in $\mathcal{O}_d$ is denoted $\CS(\bR^d)$.
\end{defn}

Every constructible set is triangulable \cite{van1998tame}, which allows us to define algebraic topological invariants such as Euler characteristic and homology. Importantly, any shape that can be faithfully represented via a mesh on a computer is a constructible set. 

\subsection{Persistence Modules}
\label{sec:pers-modules-review}

For a constructible set $M \in \CS(\mathbb{R}^d)$, the homology group $H_n(M; \Bbbk)$ in degree $n$ with coefficients in a field $\Bbbk$ provides an algebraic summary of $M$. Persistent homology extends this concept by capturing the evolution of homology groups across an $\mathbb{R}$-indexed filtration of $M$. Let $\bR$ denote the category whose objects are the real numbers and which admits a unique morphism $a \to b$ whenever $a \leq b$. 

\begin{defn}[Persistence Modules]
    A persistence module is a functor $U : \bR \to \Vect$, where $\Vect$ is the category of vector spaces over a fixed field $\Bbbk$.
\end{defn}

In other words, a persistence module $U$ assigns to each real number $t$ a vector space $U(t)$, and to each pair $s \leq t$ a linear map $U(s \leq t) : U(s) \to U(t)$, such that $U(t \leq t)$ is the identity and $U(r \leq t) = U(s \leq t) \circ U(r \leq s)$ for $r \leq s \leq t$. 

Persistence modules arise naturally from filtered topological spaces. Given a constructible set $M \in \CS(\bR^n)$ and a constructible function\footnote{A constructible function is one with a constructible graph.} $h : M \to \bR$, we obtain a filtration of $M$ by considering the sublevel sets $M_t = h^{-1}(-\infty, t]$. For each homology degree $n$, the homology groups of these sublevel sets, together with the maps induced by inclusions, form a persistence module $\text{PH}_n : \bR \to \text{Vect}$ that sends $t \mapsto H_n(M_t)$.

It is a remarkable fact of representation theory \cite{crawley2015decomposition} that whenever this functor lands in the sub-category of finite-dimensional vector spaces and maps, then the persistent module decomposes into a collection of indecomposable ``building block'' functors called interval modules, i.e., $U \cong \bigoplus_I I^{n_I}_{\Bbbk}$, where each $I^{n_I}_{\Bbbk}:(\bR,\leq)\to\Vect$ assigns $\Bbbk^{n_I}$ to points in the interval $I=[b,d)$, the identity map to $I^{n_I}_{\Bbbk}(s\leq t)$ whenever $s,t\in [b,d)$, and outside of this interval the vector spaces and maps are $0$. 
Each interval $[b,d)$ represents a homological feature that appears at time $b$ and disappears by time $d$.
We remark that the interval $I$ could also have one of the three forms: open-open $(b,d)$, open-closed $(b,d]$, and closed-closed $[b,d]$, but in practice only closed-open intervals $[b,d)$ appear.

\subsection{The Interleaving Distance}\label{sec:interleaving-distance}
The category of persistence modules, 
$\Vect^\bR$, consists of functors of the form $\bR \to \Vect$ with morphisms as natural transformations. 
This category is equipped with an interleaving distance $d_{\mathrm{Int}}$, which measures the similarity between persistence modules.
We note that this distance measures ``approximate isomorphism'' between persistence modules.

\begin{defn}[Interleaving distance]\label{def:interleaving-distance}
    Given $\varepsilon \geq 0$, one defines the \emph{$\varepsilon$-shift} operation as the map $T_\varepsilon:\bR\to \bR$ where $T_\varepsilon(x)=x+\varepsilon$. This is an endofunctor of $\bR$, when $\bR$ is viewed as a category.
    Pre-composition of a persistence module $U\in \Vect^{\bR}$ with $T_{\varepsilon}$ is called the $\varepsilon$-shift of $U$.
    An \emph{$\varepsilon$-interleaving} between two persistence modules $U$ and $V$ is simply a pair of morphisms from each original module to the other one's shift, i.e., a pair of natural transformations $\Phi : U \Rightarrow VT_\varepsilon$ and $\Psi : V \Rightarrow UT_\varepsilon$ such that $(\Psi T_\varepsilon)\Phi = U\sigma_{2\varepsilon}$ and $(\Phi T_\varepsilon)\Psi = V\sigma_{2\varepsilon}$, where $\sigma_{2\varepsilon}$ is the canonical natural transformation from the identity to $T_{2\varepsilon}$. 
    The \emph{interleaving distance} between $U$ and $V$ is then defined as:
\[
d_{\text{Int}}(U, V) = \inf\{\varepsilon \geq 0 \mid U \text{ and } V \text{ are } \varepsilon\text{-interleaved}\}
\]
If there is no interleaving between $U$ and $V$, we set $d_{\mathrm{Int}}(U,V)=\infty.$
\end{defn}

Persistence modules can be generalized by considering different choices of indexing categories, see the theory of \emph{generalized persistence modules} \cite{bubenik_generalisedPM}.  
For the purposes of our paper, we are particularly interested in the case where the indexing category is an interval $ I = [a,b] \subseteq \bR$. 
More generally, for any interval $I\subseteq \R$, we have the following definition.

\begin{defn}[$I$-persistence module]
    An $I$-persistence module is a functor $U : I \to \Vect$, where $I \subseteq \bR$ is an interval of $\bR$.
\end{defn}

The category of $I$-persistence modules denoted as $\Vect^I$.
The interleaving distance between $I$-persistence modules can be defined by letting the translation functor $T_\epsilon(x)=\min(x+\epsilon,b)$ so that translations by $\epsilon$ are capped at $b$.


\subsection{Persistent Homology Transform}
\label{sec:PHT-review}

We now introduce the central object of our study, the Persistent Homology Transform (PHT). 
The PHT captures the persistent homology of a constructible subset $M \in \CS(\bR^d)$ by considering filtrations induced by height functions in all directions. 
Specifically, for each direction $v$, we define a height function $h_v : M \to \bR$ as $h_v(x) = x \cdot v$. 
This function induces a filtration of $M$ by sublevel sets 
$$M_{v,t} = \{x \in M \mid x \cdot v \le t\} = h_v^{-1}(-\infty, t].$$

\begin{defn}[PHT: Version 1]\label{defn:map version}
	Let $M\in \CS(\bR^d)$ be a constructible set. The \define{persistent homology transform} of $M$ in degree $n$ is defined as the map
	\[
	\pht^\bR_n(M): \bS^{d-1} \to \Vect^\bR \qquad v \mapsto (t\mapsto H_n(M_{v,t}) )
	\] where $ H_n$ denotes the $n$-th homology functor.
\end{defn}

\begin{rem}[Rescaling Shapes]\label{rem:rescale}
    As observed in \cite[Rmk. 2.2]{turner} and reviewed in \cref{app: pht lipschitz proof}, the PHT is Lipschitz with respect to the Euclidean metric on $\bsd$ and the interleaving distance on $\Vect^\bR$.
    However, this Lipschitz constant is $\sup_{x\in M} ||x||_2$, which is the largest distance between $0$ and any point in $M$.
    To provide a uniform Lipschitz bound across all shapes, we assume that every shape $M$ in this paper has been rescaled so that $M$ is contained in the unit disc $\mathbb{D}^d:=\{x\in \R^d\mid ||x||_2\leq 1\}$.
\end{rem}

Building off \cref{rem:rescale}, we observe that PHT: Version 1 can be rewritten as $\pht^\mathbb{I}_i(M)$, where $\mathbb{I} = [-1,1]$ contains the range of height values where topological changes can occur. 
This observation allows us to introduce an alternative coordinate system for the PHT based on the angle $\theta$ between the direction $v$ and $x\in \bsd$ such that $x\cdot v=t$. 
We adopt the convention of using $-\theta$ to represent the angular defect from direction $v$ so that $\theta=0$ corresponds to the case where $t=1$, and so $M_{v,t=1}=M$, and $t=0$ corresponds to $\theta = -\pi/2$. 
Because we want increasing values of $t$ to correspond to increasing values of $\theta$, the choice of $\Theta = [-\pi,0]$ as the principal domain for $\arccos$ is most natural here.
We solidify this second formulation of the PHT, as it will lead to improved stability constants later in the paper.

\begin{defn}[PHT: Version 2]\label{defn:map version2}
	Let $M\in \CS(\bR^d)$ be a constructible set. The \define{angular persistent homology transform} of $M$ in degree $n$ is defined as the map
	\[
	\pht^\Theta_n(M): \bS^{d-1} \to \Vect^{\Theta} \qquad v \mapsto (\theta \mapsto H_n(M_{v,\cos{\theta}}) )
	\]

\end{defn}


\begin{rem}[Change of coordinates]\label{rem: change of coords}
    For any $v\in \bsd,\theta\in [-\pi,0]$ and $t\in [-1,1]$, the angular PHT can also be written as $\pht^\Theta_n(M)(v)(\theta)=\pht^\mathbb{I}_n(M)(v)(\cos{\theta})$ and similarly the ordinary PHT can be written as  $\pht^\mathbb{I}_n(M)(v)(t)=\pht_n^\Theta(M)(v)(-\arccos{t})$. 
    Note that the negative sign comes from the the convention of taking angles in $[-\pi,0]$.
\end{rem}

\subsection{Distances on PHTs}

We can define a distance for the PHTs of two shapes $M$ and $N$ by considering the maximum interleaving distance over all directions on the sphere; this comports with various $L^{\infty}$-type stability results. 
Since we have two versions of the PHT--- \cref{defn:map version} and \cref{defn:map version2} above---we introduce some extra notation to allow us to streamline our results.
Let $X$ be an interval of $\bR$, e.g., 
$X=\bI=[-1,1]$ or $X=\Theta=[-\pi,0]$. 

\begin{defn}[PHT interleaving distance]\label{def:PHT-distance}
    Let $M,N\in \CS(\bR^d)$. The \define{PHT interleaving distance} between $\pht^X_n(M)$ and $\pht^X_n(N)$ is defined as
    $$d_I(\pht^X_n(M), \pht^X_n(N))=\max_{v\in \bsd}d_{\mathrm{Int}} \left(\pht^X_n(M)(v), \pht^X_n(N)(v)\right),$$ where $d_{\mathrm{Int}}$ denotes the interleaving distance between persistence modules. 
\end{defn}

\begin{rem}\label{rem: bottlenck}
    When $X=\bR$ we have the usual persistent homology transform, and by the isometry theorem \cite{isometrythm}, the PHT interleaving distance coincides with the PHT bottleneck distance as defined in \cite{arya2024sheaf}.
\end{rem}

We equip $\bsd$ with two metrics, which reflects our choice of $X=\bI$ or $X=\Theta$. 
The choice of $X=\bI$ leads naturally to stability results for the Euclidean metric $d_2(v,w) = \|v-w\|_2$ on the sphere.
By contrast, the choice of $X=\Theta$ leads naturally to stability results stated in terms of the round (or geodesic) metric $d_g(v,w) = \arccos(v \cdot w)$. 

\subsection{Lipshitz Stability of the PHT}

It was shown in \cite{turner} that the PHT is 1-Lipschitz (and hence also continuous) with respect to the interleaving distance. It is straightforward to see (cf.~\cref{lem: d_g lipschitz}) that the angular-coordinatized version $\pht^\Theta$ is also 1-Lipschitz. The proof of the following lemma is provided in \cref{app: pht lipschitz proof}.

 \begin{lem}\label{lem: pht lipschitz}
    $\pht_n^\bI(M)$ and $\pht_n^\Theta(M)$ are 1-Lipschitz with respect to $d_2$ and $d_g$. 
\end{lem}

\subsection{Nesting Relations between Sub-Level Sets}\label{subsec:nesting}

We note here the relationships between sublevel sets of a shape $M\subseteq \bdd$ when viewed in different directions. A key observation is that the Euclidean metric is overly conservative for witnessing inclusions of sublevel sets. By contrast, the round metric is tight.

\begin{lem}\label{lem:d_2 lipschitz}
    For directions $v,w\in \bS^{d-1}$ and heights $s, t\in [-1,1]$, suppose $t+ d_2(v,w) \le s$ then $M_{v,t}\subseteq M_{w,s}$. 
\end{lem}
\begin{proof}
    Let $x\in M$ such that $x\cdot v\le t$. Then $|x\cdot v-x\cdot w| \le \|v-w\|_2 \le s-t$, and so $x\cdot w\le s-t + x\cdot v\le s-t+t\le s$.
\end{proof}

It is important to note that the converse is not always true. 
If $M_{v,t} \subset M_{w,s}$ then $d_2(v,w)$ may be larger than $s-t$; as the reader is encouraged to find for themselves.
This is one of several reasons to consider the more natural round (geodesic) metric on the sphere, where a converse result can be proved; cf.~\cref{lem:d_g lipschitz converse}.

\begin{lem}\label{lem: d_g lipschitz}
     For directions $v,w\in \bS^{d-1}$ and angles $\theta, \phi \in [-\pi,0]$ where $\phi+ d_{g}(v,w) \le \theta$, then $M_{v,\cos{\phi}}\subseteq M_{w,\cos{\theta}}$. 
\end{lem}
\begin{proof}
    Let $x\in M_{v,\phi}$ i.e., $x\cdot v\le \cos \phi$. Taking arc cosine on both sides gives,
    $\arccos{x\cdot v} \ge \arccos{\cos{\phi}} = -\phi$ since $\phi\in [-\pi,0]$. We are further given that $d_g(v,w) \le \theta-\phi$ and so $\arccos{x\cdot v}\ge -\phi \ge \arccos{v\cdot w}-\theta$. Rearranging gives $\theta \ge \arccos{v\cdot w}- \arccos{x\cdot v} \ge -\arccos{x\cdot w}$ by \cref{lem:triangle ineq for geodesic}. Taking cosine on both sides yields $x\cdot w\le \cos \theta$. 
\end{proof}


\begin{lem}\label{lem:d_g lipschitz converse}
    Let $M=\mathbb{D}^d$, which is the closed unit ball. 
    Given arbitrary  $v,w,\theta$, the smallest $\phi$ (in absolute value) satisfying $M_{v,\cos\phi} \subseteq M_{w,\cos\theta}$ is $\phi = \theta - \arccos(v \cdot w)$.
\end{lem}


\section{Lipschitz Extensions of the PHT}\label{sec: lipschitz extension}


We now proceed to answer our first question from the introduction: How do we take a finite collection of persistence modules $U_1,\ldots, U_n$, each of which we assume records the persistent homology of some unknown shape $M\subseteq \mathbb{D}^d$ when filtered along directions $v_1,\ldots, v_n\in \bsd$, e.g.,  $U_i\cong H_n(M_{v_i,t})$, and predict the persistent homology for an unknown direction $w\in \bS^{d-1}$?
The method we advance here is to adapt the work of \cite{de2013geometry,bubenik} on Kan extensions to the specifics of the PHT.
To contextualize the significance of our results, we first highlight the much more general result of \cite{de2013geometry}, which appeared more than 10 years ago.

\begin{thm}[Thm. 4.1 of \cite{de2013geometry}]\label{old deSilva theorem}
    Suppose $A$ is a subspace of a metric space $(Z,d_z)$.
    Any \textbf{coherent map} $f:A\to \Vect^{\R}$ can be extended to a 1-Lipschitz map $\hat{f}:Z \to \Vect^{\R}$.
\end{thm}

The definition of a coherent map\footnote{As introduced in \cite{de2013geometry}, the definition of coherent consists of two properties: (1) $f$ must be 1-Lipschitz and (2) there must be a system of interleavings between $f(a)$ and $f(b)$ for all pairs $a,b\in A$, which further satisfies a condition that must be checked on all triples $a,b,c\in A$. As we'll see, (2) implies (1).} is deferred to \cref{lem: 1lipschitz category}, where a cleaner categorical interpretation is used.
For now, we stress that the work of \cite{de2013geometry} pre-dated the introduction of the persistent homology transform in \cite{turner}, so it was not in the scope of possible applications back then. 
Although the more robust follow-up paper \cite{bubenik} was published several years after \cite{turner}, our work here seems to be the first place where the methods of \cite{de2013geometry,bubenik} are applied to the study of the PHT.

Moreover, on a technical side, we note that we have to modify \cref{old deSilva theorem} and much of \cite{bubenik} to work with $X$-modules where $X=\bI$ or $X=\Theta$.
Additionally, we observe in \cref{prop:PHT-coherence} that the coherence condition is automatically satisfied in the PHT setting, 
i.e., every map $f:A\to \Vect^X$, where $f(a)$ is the PHT in direction $a\in A$, is automatically coherent and thus admits a Lipschitz extension.
To better understand this extension, we emphasize that extra data on how the modules $f(a)$ and $f(b)$ interact is required.
This ``interaction'' between observed persistence modules, already hinted at in \cref{subsec:nesting}, is best captured by the introduction of an auxiliary space-time category, which was one of the main theoretical contributions that separated \cite{bubenik} from \cite{de2013geometry}.
The actual extension process is then accomplished via an elegant application of Kan extensions reviewed below and in \cref{app:Kan}.

\subsection{The Space-Time Category}

We now proceed with adapting the ideas from \cite{bubenik} to our setting.
Let $d$ be a distance metric on $\bsd$. 
Let $X\subseteq \bR$ be an interval of $\bR$ such that $\sup X\in X$ and let $A$ be a subset of $\bsd$.
We define an auxiliary pre-order category\footnote{We use the notation $\AX$ to emphasize that this object is typed as a category. Often we will need to pass back and forth from maps from $A\subseteq \bsd$, viewed as a set, and functors from $\AX$.} $\AX$, whose objects are elements of the set $A\times X$.
We declare the existence of a unique morphism $(v,\alpha) \to (w,\beta)$ whenever $\alpha+ d(v,w)\le \beta$. 
Additionally, for any \( \alpha \in X \) and for any \( v, w \in S^{d-1} \), there is always a unique morphism \( (v, \alpha) \to (w, \sup X) \), where \( \sup X \) is the maximum value in \( X \). 
Geometrically, this additional morphism is motivated by the observation that \( M_{v, \alpha} \subseteq M_{w, \sup X} \) always, since \( M_{w, \sup X} =M \) is the whole shape.
Algebraically, these additional morphisms imply that $\AX$ is a pre-order and not a poset, because for arbitrary directions $v,w\in A$ we will have a morphism $(v,\sup X)\to (w,\sup X)$ and vice versa, but $v$ and $w$ need not be the same.
As usual, we consider two specializations of the metric $d$ and the interval $X\subset \bR$:

\medskip

\underline{\textbf{The Euclidean case, where $d = d_2$ and $X = \mathbb{I}=[-1,1]$:}}

In this case, the category $\bA\bI$ consists of objects $(v,t)$ with the pre-order relation $(v,t) \leq (w,s)$ if and only if $t+ d_2(v,w) \leq s$. 
Geometrically, by \cref{lem:d_2 lipschitz}, whenever there is an arrow $(v,t) \leq (w,s)$ in the category, we have an inclusion of spaces $M_{v,t} \subseteq M_{w,s}$. Additionally, for any $t\in \bI$ and $v,w\in \bsd$, there is always a morphism from $(v,t)\le (w,1)$. 

\medskip

\underline{\textbf{The Geodesic case, where $d = d_g$ and $X =\Theta=[-\pi,0]$:}}
    
Here, $\bA\Theta$ consists of objects $(v,\phi)$ with the pre-order relation $(v,\phi) \leq (w,\theta)$ if and only if $\phi+ d_g(v,w) \leq \theta$. 
By \cref{lem: d_g lipschitz}, whenever there is a morphism $(v,\phi) \leq (w,\theta)$, we have an inclusion of spaces $M_{v,\cos\phi} \subseteq M_{w,\cos\theta}$. Additionally, for any \( \theta \in [-\pi, 0] \) and any \( v, w \in \bS^{d-1} \), we formally add the morphism 
$(v, \theta) \leq (w, 0)$, since \( M_{v, \cos\theta} \subseteq M_{w, \cos{0}} = M \).

\medskip
    


As noted in \cite{bubenik}, the coherent definition of \cite{de2013geometry} can be more succinctly stated by requiring that a map $f:A\to\Vect^X$ arise as a restriction of a functor $G:\AX\to \Vect$.
Such maps $f$ are then automatically 1-Lipschitz.
This also provides a novel proof of Lipschitz stability of the PHT, cf. \cref{cor:lipschitz-PHT}.

\begin{lem}[cf.~\cite{bubenik}]
\label{lem: 1lipschitz category}
    Let $f: A\to \Vect^X$ be a map from any metric space $(A,d_A)$, where $X\subseteq \bR$ is an interval that includes its supremum. 
    If there is a functor $G: \AX \to \Vect$ such that $f(v)=G(v,\cdot)$ for all $v\in A$, then $f$ is said to be \textbf{coherent} and is automatically $1$-Lipschitz with respect to the interleaving distance on $\Vect^X$. 
\end{lem}

\begin{proof}
The definition of coherence, adapted from \cite{bubenik}, is simply that $f(v)=G(v,\cdot)$ for some functor $G:\AX\to\Vect$.
We now prove that coherent maps are $1$-Lipschitz with respect to the modified interleaving distance for $X$-modules.
    Let $v,w\in A$ be any two distinct points and let $\epsilon=d_A(v,w)$.
    Recall the translation functor on $X$ is defined by $T_\epsilon(x)=\min(x+\epsilon,\sup X)$. 
    We need to prove that an $\epsilon$-interleaving exists between $f(v)$ and $f(w)$.
    To this end we must show that for any $\alpha \in X$ there are morphisms in $\AX$ 
    $$(v,\alpha)\le (w,T_\epsilon(\alpha))\qquad \mathrm{and}\qquad (w,T_\epsilon(\alpha))\le (v,T_{2\epsilon}(\alpha)).$$
    We break this up into two cases.
    When $\alpha+2\epsilon< \sup X$, these morphisms exist by the definition of $\AX$. 
    For the case where $\alpha+2\epsilon \ge  \sup X$, we have defined $\AX$ in such a way so that $(v,\alpha) \leq (w,\sup X)$ always exists. 
    Since $\AX$ is a thin category, all diagrams commute. Applying the functor $G$ to these morphisms yields the desired $\epsilon$-interleaving between $G(v,\cdot)$ and $G(w,\cdot)$. 
    Since $f(v)=G(v,\cdot)$ and $f(w)=G(w,\cdot)$, we have that $d_I(f(v),f(w)) \le \epsilon=d_A (v,w)$. Since $v$ and $w$ were arbitrary, this proves the result. 
\end{proof}


\begin{prop}[Coherence of the PHT]\label{prop:PHT-coherence}
    If $A\subseteq \bsd$ is a subset of directions and $f:A\to\Vect^X$ is the persistent homology transform of a shape $M\subseteq \bdd$, when viewed in either Euclidean or angular coordinates, then $f$ is automatically coherent.
\end{prop}
\begin{proof}
    We prove the result for $A=\bsd$ since any $\AX$ can be viewed as a restriction of $\bsd\bX$.
    For $\mathrm{PHT}_n^\bI(M): \bsd \to \Vect^\bI$, we can define $G^\bI: \bsd\bI \to \Vect$ by $G^\bI(v,t) = H_n(M_{v,t})$. Similarly, for $\mathrm{PHT}_n^\Theta(M): \bsd \to \mathrm{Vect}^{[-\pi,0]}$, we define $G^\Theta: \bsd[-\pi,0] \to \mathrm{Vect}$ by $G^\Theta(v,\phi) = H_n(M_{v,\cos\phi})$. In both cases, one can verify that $G^\bI$ and $G^\Theta$ are functors on $\bsd\bI$ and $\bsd\Theta$ using \cref{lem:d_2 lipschitz} and \cref{lem: d_g lipschitz}, respectively.
\end{proof}

\begin{cor}\label{cor:lipschitz-PHT}
     \cref{lem: pht lipschitz} is true, i.e. both $\PHT_n^{\bI}(M)$ and $\PHT_n^{\Theta}(M)$ are 1-Lipshitz.
\end{cor}

    

\subsection{Left Kan Extensions for the PHT}

\cref{prop:PHT-coherence} and \cref{old deSilva theorem} together imply an existence result: Given an \emph{arbitrary} subset $A\subseteq \bsd$ one can extend the PHT at $A$ to the entire sphere.
This does not explain \emph{how} to compute such an extension, however.
In this section, we show how the general technology of Kan extensions (see \cref{app:Kan}) provides a constructive solution to this problem.
As noted in \cite{bubenik}, there are three possible choices of Kan extensions: left, right, and center. 
We will focus on the left Kan extension in this section and a simple example will motivate consideration of the other two.

The general setup goes as follows: let $G^X$ be a functor from $\AX$ to $\mathrm{Vect}$. 
We seek an extension $\hat{G}^X$ as illustrated in the diagram below:
\[
\begin{tikzcd}
	\AX & \mathrm{Vect} \\
	\mathbb{S}^{d-1}\mathbb{X}
	\arrow[hook, from=1-1, to=2-1]
	\arrow["G^X", from=1-1, to=1-2]
	\arrow["{\hat{G}}^X"', dashed, from=2-1, to=1-2]
\end{tikzcd}
\]

When $X = \mathbb{I}$ and $d=d_2$, we define $G^\mathbb{I} = \mathrm{PHT}_n^\bI(M)$, and when $X = \Theta$ and $d=d_g$, we define $G^\Theta(v,\phi) = \mathrm{PHT}_n^\Theta(M)$. 
Either choice of coordinate system allows us to frame the extension of the PHT as a Kan extension problem. 
Regardless of one's choice, we write:
\[
\overline{\pht}_n^X(M):\bsd\to \Vect^X\qquad v\mapsto \widehat{\pht}^X_n(M)(v,\cdot)
\]

\begin{figure}
    \centering
    \includegraphics[width=0.65\linewidth]{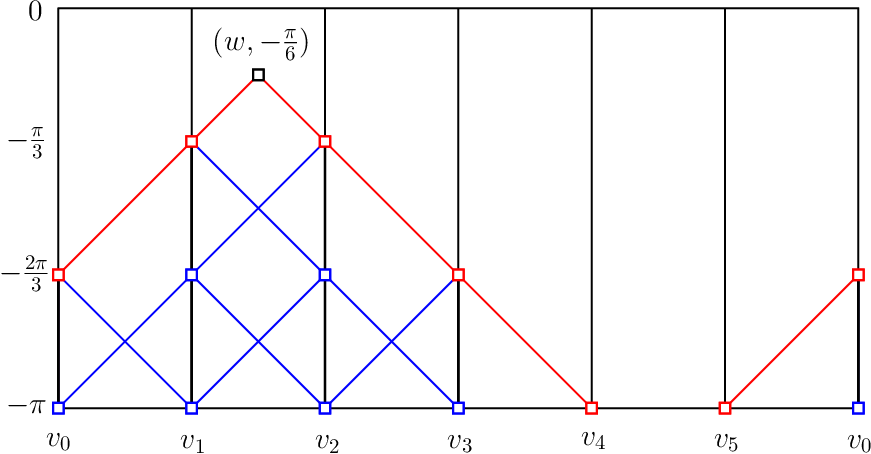}
    \caption{Given directions $A=\{v_k=e^{2\pi i k/6}\mid k=0,\ldots,5\} \subset \bS^1$, marked on the horizontal axis, our goal is to calculate the value of the PHT at the unknown direction $w=e^{i\pi/2}$ and $\Theta$-value $-\pi/6$. To do this, we look at the past light cone, whose boundary is drawn in red. The red boxes indicate the largest $\Theta$-values for each direction that can include into $M_{w,\cos(-\pi/6)}$. The blue lines indicate the past light cones of each of these maximal values, which will be important when identifying common topological features between directions.}
    \label{fig:past-light-cone}
\end{figure}

For an object $(w,\beta)\in \bS^{d-1}\bX$, the left Kan extension examines the \define{past light cone}:\footnote{For the reader familiar with special relativity, the past light cone encodes all points in space-time that can ``causally'' influence an event at $(w,\beta)$ by ensuring that each $(v,\alpha)\in L^-(w,\beta)$ satisfies $d(v,w)\le \beta - \alpha$. 
The Lipschitz-1 condition is exactly the restriction that nothing moves faster than the speed of light $c=1$.}
\[
L^-(w,\beta) \vcentcolon= \{ (v,\alpha) \in \AX \mid \alpha + d(v,w) \le \beta \}
\] 
This is simply the collection of sub-level sets in our ``known'' subcategory $\AX$ that includes into the sub-level set in direction $w$ and height $\beta$; see \cref{fig:past-light-cone}.
More formally, $L^-(w,\beta)$ is the projection of the under comma category $\iota \downarrow (w,\beta)$ to $\AX$ where $\iota \colon \AX \hookrightarrow \bS^{d-1}\bX$ is the obvious inclusion.
Following \cref{thm:ptwise-Kan}, the left Kan extension of $G^X$ along $\iota$ can be computed via the colimit
\[
\widehat{\PHT}_n^{X}(M)(w,\beta) \vcentcolon= \operatorname*{\underrightarrow{\mathrm{colim}}}_{(v,\alpha) \in L^-(w,\beta)} \PHT_n^X(M)(v,\alpha).
\]
As we will see in the following examples, this colimit can be thought of as the ``best guess'' using the homologies of known sub-level sets and any relations between these.

\begin{figure}
    \centering
    \includegraphics[width=0.65\linewidth]{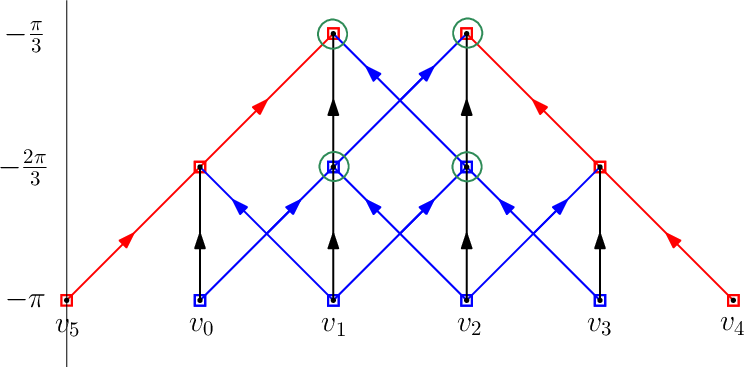}
    \caption{Continuing the example from \cref{fig:past-light-cone}, we see that the two maximal sub-level sets at $(v_1,-\pi/3)$ and $(v_2,-\pi/3)$ are related via blue arrows. In the case considered in \cref{ex:disk_geod_left}, of predicting the PHT in degree 0 for the shape $M=\bD^2$, this diagram reduces to constant connected diagram whose value everywhere is $\Bbbk$. The colimit will also be $\Bbbk$.}
    \label{fig:colim-light-cone}
\end{figure}

\begin{exmp}[Left Kan Extension, Correct Value]\label{ex:disk_geod_left}
    For the simplest non-trivial example, consider the problem of interpolating the PHT for the shape $M=\bD^2$ using only the known directions $A=\{v_k=e^{2\pi i k/6}\mid k=0,\ldots,5\} \subset \bS^1$.
    As $M$ is a simply connected subset of the plane, only $H_0$ is of interest here.
    As initiated in \cref{fig:past-light-cone}, we are interested in computing the predicted homology in degree 0 at the sublevel set $(w=e^{i\pi/2},\beta=-\pi/6)$.
    In this case, the past light cone $L^-(w,\beta)$ has two maximal sub-level sets at $(v_1,-\pi/3)$ and $(v_2,-\pi/3)$.
    A standard recipe for computing the colimit of $H_0(M_{-,-})$ over $L^-(w,\beta)$ is to identify these two maximal sub-level sets as part of a cofinal subcategory. 
    Following \cref{defn:cofinal}, one can verify that the top two red boxes \emph{along with} the top two blue boxes form a four-object co-final subcategory $\cC(w,\beta)$ of $L^-(w,\beta)$; circled in green in \cref{fig:colim-light-cone}.
    Specifically, $\cC(w,\beta)$ consists of the blue crossed arrows $(v_2,-2\pi/3)\to (v_1,-\pi/3)$ and $(v_1,-2\pi/3)\to (v_2,-\pi/3)$ along with the black internal arrows $(v_1,-2\pi/3)\to (v_1,-\pi/3)$ and $(v_2,-2\pi/3)\to (v_2,-\pi/3)$.
    Since $H_0(M_{v_1,\cos -\pi/3})=\Bbbk=H_0(M_{v_1,\cos -2\pi/3})$ and $H_0(M_{v_2,\cos -\pi/3})=\Bbbk=H_0(M_{v_2,\cos -2\pi/3})$ the colimit of this four object diagram is equivalently the cokernel of the map
    \[
    W: H_0(M_{v_1,\cos -2\pi/3})\oplus H_0(M_{v_2,\cos -2\pi/3}) \to H_0(M_{v_1,\cos -\pi/3})\oplus H_0(M_{v_2,\cos -\pi/3}),
    \]
    which in the standard basis is
    \[
    \begin{bmatrix}
        1 & 1 \\
        1 & 1
    \end{bmatrix}.
    \]
    The cokernel, and the predicted value for $H_0(M_{w,\cos(-\pi/6)})$ is $\Bbbk$, which is correct.
\end{exmp}

\begin{exmp}[Left Kan Extension, Wrong Value]\label{ex:left-kan-wrong-answer}
    Continuing \cref{ex:disk_geod_left}, consider the same direction $w=e^{i\pi/2}\in \bS^1$ as before, but with $\beta=-2\pi/3$.
    Following \cref{fig:disconnected-light-cone}, the reader should convince themselves that $L^-(w,\beta)\cong [(v_1,-\pi),(v_1,-5\pi/6)]\sqcup [(v_2,-\pi),(v_2,-5\pi/6)]$ is the disjoint union of two categories, each of which are intervals.
    This contains a co-final subcategory $\mathcal{C}(w,\beta)\subseteq L^-(w,\beta)$ consisting of the two objects $(v_j,-5\pi/6)$ for $j=1,2$.
    As such, the left Kan extension can be computed as
    \[
    \overline{\pht}_0^X(M)(w)(\beta)\cong H_0(M_{v_1,\cos(-5\pi/6)})\oplus H_0(M_{v_2,\cos(-5\pi/6)})\cong \Bbbk\oplus \Bbbk.
    \]
    This disagrees with the correct answer of $\Bbbk$, which is the homology of the disk intersected with an overlapping half-plane.

\begin{figure}
    \centering
\includegraphics[width=0.65\linewidth]{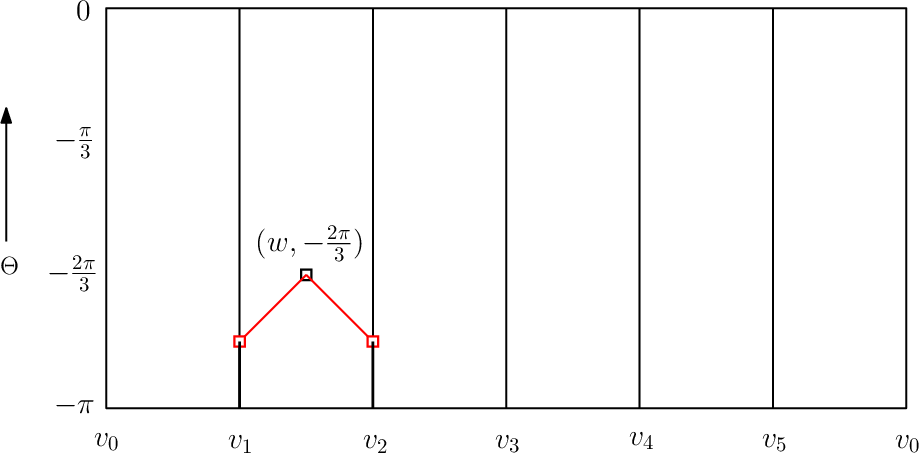}
    \caption{The past light cone at the value $w=e^{i\pi/2}$ and $\beta=-2\pi/3$ is that part of $\bA\Theta$ that is below the red cone. This corresponds to the two intervals $[(v_j,-\pi),(v_j,-5\pi/6)]$ for $j=1,2$. Because $d_g(v_1,v_2)=\pi/3$ there can be no internal morphisms between these sub-categories, so the past light cone is disconnected. This will cause the Kan extension of any constant functor to be non-constant.}
    \label{fig:disconnected-light-cone}
\end{figure}
\end{exmp}

One can see from \cref{ex:disk_geod_left} and \cref{ex:left-kan-wrong-answer} that the left Kan extension tends to produce spurious features for small values of $\beta$. 
This is primarily due to topological aspects of the past light cone, which we explain in the following remark.

\begin{rem}[Empty and Disconnected Components from Left Kan Extensions]\label{rem:disconnected-light-cone}
    Suppose $A=\{v_1,\ldots, v_n\}\subset \bsd$ is a finite set of directions.
    The past light cone $L^-(w,\beta)$, equivalently the under comma category of $\iota:\AX \hookrightarrow \bsd\bX$ at $(w,\beta)$, will only include information from the sub-level sets for directions $v_i$ that satisfy
    $d(w,v_i)\leq \beta-\inf X$.
    Notice that as $\beta\to \inf X$ this will eliminate more and more directions to the point that when $\beta-\inf X$ is less than $\min\{d(v_i,w)\}$, the past light cone will be empty.
    For such a value of $\beta$, the left Kan extension will assign the zero vector space.

    Now suppose we are in a range of $\beta$ where there are only two directions $v_1$ and $v_2$ such that $d(v_1,w)\leq d(v_2,w)\leq \beta-\inf X$.
    At these values of $\beta$, the left Kan extension will assign the homology of the disjoint union (i.e., the sum of the homologies) unless $\beta$ also satisfies
    \[
    d(v_1,v_2)+d(v_2,w) \leq \beta - \inf X.
    \]
    The gap between $d(v_2,w)$ and $d(v_1,v_2)+d(v_2,w)$ indicates a region where the past light cone is disconnected.
    It should be noted that the extra factor $d(v_1,v_2)$ is the length of the shortest morphism needed to formally ``glue'' homological features observed in direction $v_1$ with features observed in direction $v_2$.
\end{rem}

\subsection{Right and Center Kan Extensions for the PHT}

As \cref{rem:disconnected-light-cone} illustrates, most directions $w$ and small values of $\beta$ will produce an empty or disconnected past light cone.
This leads to undesirable estimates for the homology in those parameter ranges---either zero or a direct sum of homologies, respectively.
The reader who is sensitive to duality may ask if there is a way to exploit the \define{future light cone}
\[
L^+(w,\beta) \vcentcolon= \{ (v,\alpha) \in \AX \mid \beta+ d(v,w) \le \alpha \},
\]
which tends to be non-empty and connected for values of $(w,\beta)$ where $L^-(w,\beta)$ is empty or disconnected.
The \emph{right} Kan extension follows a similar process as the left Kan extension, except one has to exchange past with future and colimits for limits.
Specifically, the right Kan extension assigns to $(w,\beta)$ the value
\[
\varprojlim_{(v,\alpha) \in L^+(w,\beta)} \PHT^X(M)(v,\alpha).
\]
This assignment is slightly counterintuitive, because formally to construct this limit, one would need to take the product of all homologies of ``future'' sub-level sets that $(w,\beta)$ might map to, and then assign a best approximating sub-vector space of this product.

\begin{figure}
    \centering
    \includegraphics[width=0.5\linewidth]{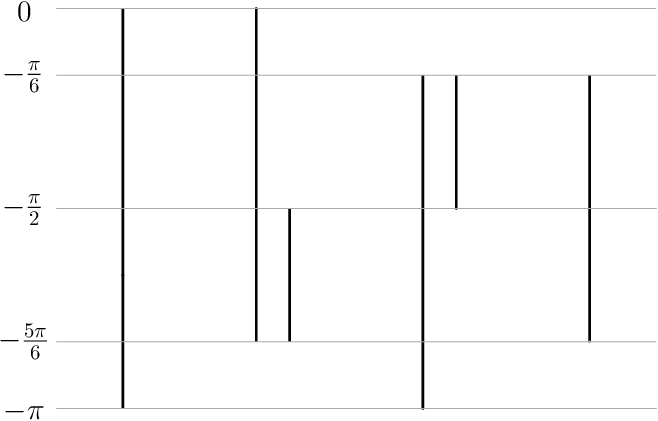}
    \caption{From left to right: 
    (1) the true $\PHT^\Theta_0$ at direction $w=e^{i\pi/2}\in \bS^1$,
    (2) the left-Kan extension of $\PHT^\Theta_0$ at \(w\) using $A=\{v_k=e^{2\pi i k/6}\mid k=0,\ldots,5\}$,
    (3) the right-Kan extension for the same $w$ and $A$, finally (4) the center-Kan extension. See \cref{ex:Kan-disc} for an explanation.}
    \label{fig:barcode_geo}
\end{figure}

\begin{exmp}[Right Kan Extension for $M=\bD^2$]\label{ex:right-Kan-disc}
    Continuing \cref{ex:disk_geod_left} and \cref{ex:left-kan-wrong-answer}, we consider the right Kan extension of $\PHT^\Theta_0$ for the shape $M=\bD^2$ at direction $w=e^{i\pi/2}\in \bS^1$.
    Since the shape $M=\bD^2$ is entirely symmetric, one can expect that the limit over the future light cone to be entirely dual to the colimit over the past light cone. As one can see in \cref{fig:barcode_geo}, the predicted barcode for the right Kan extension is exactly the mirror image of the predicted barcode or the left Kan extension.
\end{exmp}

As \cref{ex:right-Kan-disc} illustrates, the right Kan extension produces the correct estimate for the homology of a shape for small values of $\beta$, but has its own problems for large values of $\beta$, which is where the left Kan extension succeeds.
The \define{center Kan extension} mediates these two extremes by looking at the image of the mapping from the left Kan extension to the right.
To see why such a mapping exists, note that the limit over the future light cone (i.e., the right Kan extension) stands as a categorical co-cone to the past light cone, and thus the universal property of the colimit (i.e, the left Kan extension) ensures the existence of a unique mapping
\[
    \operatorname*{\underrightarrow{\mathrm{colim}}}_{(v,\alpha) \in L^-(w,\beta)} \PHT^X(M)(v,\alpha) \to \varprojlim_{(v,\alpha) \in L^+(w,\beta)} \PHT^X(M)(v,\alpha).
\]
The center Kan extension assigns to the object $(w,\beta)$ the \emph{image} of this natural map. 
As we will see in our example, the center Kan extension correctly ignores spurious features that both the right and left Kan extensions tend to introduce.

\begin{exmp}[Left, Right, and Center Kan Extensions for $M=\bD^2$]\label{ex:Kan-disc}
    Continuing our previous examples and examining \cref{fig:barcode_geo} we see that the left Kan extension estimates the persistent homology in direction $w=e^{i\pi/2}$ to be the sum of the interval modules $[-\frac{5\pi}{6},0]$ and $[-\frac{5\pi}{6},-\frac{\pi}{2})$.
    The right Kan extension estimates the persistent homology to be the sum of the interval modules $[-\pi,-\frac{\pi}{6}]$ and $(-\frac{\pi}{2}, -\frac{\pi}{6}]$.
    The discussion above guarantees a mapping between these modules, but the only non-zero mapping that can occur is between $[-\frac{5\pi}{6},0]$ and $[-\pi,-\frac{\pi}{6}]$.
    The image of this mapping is exactly the intersection, which is the interval $[-\frac{5\pi}{6},-\frac{\pi}{6}]$. We note that the center Kan extension is the only one whose estimates agree with an actual shape $\widehat{M}$: the regular hexagon with normal directions at $e^{(2k+1)i\pi/6}$ for $k=0,\ldots, 5$.
\end{exmp}

\subsection{Lipschitz Extensions of the PHT}
Having introduced three canonical methods for computing Kan extensions of the PHT, we record for posterity the following corollary of \cref{lem: 1lipschitz category} and \cref{prop:PHT-coherence}:

\begin{thm}[Lipschitz Extension Theorem]\label{thm:PHT-LET}
    Let $A\subseteq \bsd$ be an arbitrary subset of directions. Suppose we are given the data of a functor $G^X:\AX\to \Vect$. 
    Specifically, we consider the case where $G^X$ records the persistent homology in degree-$n$ of a set $M\subseteq \bdd$ for every direction $v\in A$, but where we are also given the maps on homology $H_n(M_{v,\alpha}) \to H_n(M_{w,\beta})$ whenever $(v,\alpha)\to(w,\beta)$ is a morphism in $\AX$.
    Any of the three extensions described aboove will produce a $1$-Lipschitz map ${\overline{\pht}}^X_n(M):\bsd\to \Vect^X$.
\end{thm}

\begin{rem}[Avoiding Maps Between Directions]\label{rem:no-maps}
    \cref{thm:PHT-LET} suffers from an apparent defect: one needs to know the maps on homology between directions in order to compute any of the three extensions described. 
    However, in practice, one is not given this information.
    \cite[\S 3.4 and \S3.5]{de2013geometry} points out that one can computationally estimate such mappings by coming up with \emph{coherent matchings} between observed persistence diagrams.
    Pairwise matchings are equivalent to maps between modules, thanks to the famed Isometry Theorem of persistence \cite{lesnick2015theory}, although this association is not always functorial \cite{bauer2014induced}.
    Using coherent matchings to compute Kan extensions is left to future work.
\end{rem}

\section{Error Bounds for Kan Approximations}\label{sec: error bounds}

In this section, we establish precise error bounds between the true PHT of a shape $M$ and the Kan extension of the PHT computed at a finite set of directions. 
This result is proved in two steps: \cref{lem: interleaving Kan} proves a point-wise bound for any direction $w\in\bsd$ and \cref{thm:stability_kan_ext} then converts this result into a distance bound between PHTs, cf.~\cref{def:PHT-distance}.
Subtle differences between the geodesic $X=\Theta$ and Euclidean $X=\bI$ coordinate systems for the PHT are exhibited through a few computations and remarks.
Finally, we generalize our error bound further by considering both a finite set of directions and scalar heights.
The machinery of Kan extensions works flawlessly here and produces a linear error term in our approximation.


\subsection{Direction-wise Bound}

\begin{lem}\label{lem: interleaving Kan}
Let $M$ be a constructible set in $\bR^d$.
For any direction $w \in \bsd$ and homological degree $n$, there exists an $\epsilon=2d(v,w)$-interleaving between $\mathrm{PHT}^X_n(M)(w)$ and $\overline{\mathrm{PHT}}^X_n(M)(w)$, where $v\in A$ is a direction where $\PHT(M)$ is sampled.
\end{lem}

\begin{rem}[Simplifying Assumptions and Notation]\label{rem:simplify-interleaving-proof}
    For the proof of \cref{lem: interleaving Kan} we will make several simplifying assumptions.
    First, we will assume that $\overline{\mathrm{PHT}}^X_n(w)$ is given by the left Kan extension, as all we really need for the proof is a construction that is functorially constructed in terms of limits and colimits.
    Additionally, we will assume $X=\bI$ since $\mathrm{PHT}^X_n(M)(w)(\beta)=H_n(M_{w,\beta})$ is more straightforward to write. Recall that when $X=\Theta$ the ``sub-level set at parameter $\beta$'' is actually $M_{w,\cos\beta}$ in the usual Euclidean coordinate system.
    Finally, as $M$ and $i\geq 0$ are constant throughout, we will ignore these as well.
\end{rem}

\begin{proof}
In view of \cref{rem:simplify-interleaving-proof} we will write $C_w:=\overline{\mathrm{PHT}}^X_n(w)$ to emphasize that $C_w(\beta)$ is given by the colimit of $\pht^X_i$ over the past light cone $L^-(w,\beta)$.
Our task is to compare this with the usual persistent homology module of $M$ in direction $w$, which we will write as $H_w$.
The reader is urged to recall \cref{def:interleaving-distance} and the definition of our translation functor $T_{\epsilon}:X\to X$, given by $T_{\epsilon}(\beta)=\min\{\beta+\epsilon,\sup X\}$.
Our task in this proof is to construct two natural transformations
\[
    \Phi: C_w \Rightarrow H_wT_{\epsilon} \qquad \text{and} \qquad \Psi: H_w \Rightarrow C_wT_{\epsilon}
\]
such that $(\Psi T_\epsilon)\Phi = C_w\sigma_{2\epsilon}$ and $(\Phi T_\epsilon)\Psi = H_w\sigma_{2\epsilon}$, where $\sigma_{2\epsilon}$ is the canonical natural transformation from the identity to $T_{2\epsilon}$. 

We first note that $\Phi:C_w \Rightarrow H_wT_{\epsilon}$ exists for every $\epsilon\geq 0$.
Indeed, since $C_w(\beta)$ is the colimit over all homology vector spaces $H_n(M_{v,\alpha})$ such that $(v,\alpha)\to (w,\beta)$ exists in the under comma category $\iota \downarrow (w,\beta)$ for the inclusion $\iota:\AX\hookrightarrow \bsd\bX$, then we can guarantee that $M_{v,\alpha}\subseteq M_{w,\beta}$ is a valid geometric inclusion, hence $H_n(M_{w,\beta})=:H_w(\beta)$ stands as a co-cone to the diagram over $L^-(w,\beta)$.
By definition, $C_w(\beta)$ is the universal co-cone over $L^-(w,\beta)$, so there is a unique map $C_w(\beta)\to H_w(\beta)$.
As this argument is functorial in $\beta$, we conclude the existence of a natural transformation $\Phi:C_w\Rightarrow H_w$ and hence a natural transformation to any shift $H_wT_{\epsilon}$ because $H_wT_{\epsilon}(\beta):=H_n(M_{w,\beta+\epsilon})$ is also a co-cone.

The construction of $\Psi:H_w\Rightarrow C_wT_{\epsilon}$ is the main difficulty of the proof, and illustrates how $\epsilon$ depends on the choice of $v\in A$.
First we note that for any pair of vectors $v,w\in \bsd$ there is always a natural transformation $\Delta_{w,v}:H_w\Rightarrow H_v T_{d(v,w)}$ because of the nesting relationship $M_{w,\beta}\subseteq M_{v,\beta+d(v,w)}$ proved in \cref{subsec:nesting}.
Moreover, if $v\in A$ there is always a natural transformation $\Sigma_{v,w}: H_v \Rightarrow C_wT_{d(v,w)}$ because $H_n(M_{v,\beta})$ participates in the diagram over $L^-(w,\beta+d(v,w))$ that defines $C_w(\beta+d(v,w))$.
Shifting this natural transformation yields a natural transformation $\Sigma_{v,w}T_{d(v,w)}:H_vT_{d(v,w)} \Rightarrow C_wT_{d(v,w)}T_{d(v,w)}=C_wT_{2d(v,w)}$.
Setting 
\[
    \Psi:=(\Sigma_{v,w}T_{d(v,w)})\Delta_{w,v}:H_w \Rightarrow H_v T_{d(v,w)} \Rightarrow C_wT_{2d(v,w)}
\]
reveals an $\epsilon$-interleaving for all $\epsilon\geq 2d(v,w)$.
\end{proof}

\begin{exmp}\label{ex:disc-direction-bound}
    Examining our running example from \cref{fig:barcode_geo}, we observe that the true $\Theta$-persistence module in direction $w=e^{i\pi/2}$ and any of the three Kan extensions are distance $\frac{\pi}{6}$ apart.
    This is less than the $\frac{2\pi}{6}$-interleaving predicted by \cref{lem: interleaving Kan}, which is the worst case scenario. 
\end{exmp}

\begin{exmp}[Euclidean Coordinates]\label{ex: disk_eu}
    If we repeat the example started in \cref{ex:disk_geod_left}, but in Euclidean $X=\bI$ coordinates, we see that \cref{lem: interleaving Kan} will predict a $2\sqrt{2-\sqrt{3}}$-interleaving, as $\sqrt{2-\sqrt{3}}$ is the Euclidean distance between $w=e^{i\pi/2}$ and $v=e^{i\pi/3}$.
\end{exmp}


\begin{exmp}[A Tight Bound]
    To see that the bound in \cref{lem: interleaving Kan} is tight, consider the shape $M=\{p=(0,1)\}$, consisting of a single point, with a single sample direction $v=e^{i0}$, pointing to the east. When querying at the south pole $w=e^{-i\pi/2}$, the true $\Theta$-persistence module for $H_0$ is born at $\theta = -\pi$, while the left-Kan extension is born at $\theta=0$. Indeed, the $\Theta$-persistence module at $v$ is born at $\theta=-\pi/2$, and only enters the past light cone of $w,\beta$ once the query angle $\beta$ satisfies $\beta \ge -\pi/2 + d_g(v,w) = 0$.  The resulting interleaving distance is therefore $\pi$. Since $d_g(v,w)=\pi/2$, this error matches the predicted worst-case bound $2(\pi/2)=\pi$. 
\end{exmp}

\subsection{Global Distance Bound}

\begin{thm}[Stability of the Kan Extension]\label{thm:stability_kan_ext}
Let $M$ be a fixed constructible set in $\mathbb{R}^d$, and let $A \subseteq \bS^{d-1}$ be a subset of directions. 
The PHT interleaving distance between the persistent homology transform $\mathrm{PHT}^X_n(M)$ and its Kan extension approximation $\overline{\mathrm{PHT}}^X_n(M)$ is bounded by twice the Hausdorff distance between $A$ and the sphere $\bS^{d-1}$ for every degree $n$, i.e.,
$$
d_I(\mathrm{PHT}^X_n(M), \overline{\mathrm{PHT}}^X_n(M)) \leq 2d_H(A, \bS^{d-1}).
$$
\end{thm}

Here, the distance metric on the sphere is the Euclidean metric when $X = \bI$, and the geodesic metric when $X = \Theta$. 
In particular, this result implies that if $A$ is an $\epsilon$-net of the sphere, then
$$
d_I(\mathrm{PHT}^X_n(M), \overline{\mathrm{PHT}}^X_n(M)) \leq 2\epsilon.
$$

\begin{proof}
Suppose $\epsilon = d_H(A, \bS^{d-1})$, which implies that for every $w \in \bS^{d-1}$, there exists a point $v \in A$ such that $d(v, w) \leq \epsilon$. By \cref{lem: interleaving Kan}, for any such pair $(v, w)$, there is a $2\epsilon$-interleaving between $\mathrm{PHT}^X_n(M)(w)$ and $\overline{\mathrm{PHT}}^X_n(M)(w)$ for every homological degree $n$. Since this holds for every $w \in \bS^{d-1}$, we conclude that the interleaving distance between the true PHT and its Kan extension approximation, measured over all directions in $\bS^{d-1}$, is at most $2\epsilon$. Thus, 
$$
d_I(\mathrm{PHT}^X_n(M), \overline{\mathrm{PHT}}^X_n(M)) = \sup_{w \in \bS^{d-1}} d_\mathrm{Int}\left(\mathrm{PHT}^X_n(M)(w), \overline{\mathrm{PHT}}^X_n(M)(w)\right) \leq 2\epsilon.
$$

\end{proof}

\begin{rem}[Bottleneck Distance Bounds]
    Following \cref{rem: bottlenck}, one could repeat this analysis for $X=\R$ and this would give a bound on the bottleneck distance for the PHT.
\end{rem}

\begin{exmp}\label{ex:}
Continuing \cref{ex:disc-direction-bound} we note that this example illustrates the bound given by \cref{thm:stability_kan_ext} because the Hausdorff distance between $A=\{v_k=e^{2\pi i k/6} \mid k=0,\ldots, 5\}$ and $\bS^1$ is half the maximal angular separation \(\frac{\pi}{3}\).
\end{exmp}

When considering global distance bounds it is worth noting that the geometry of the indexing set $X$ provides different perspectives on what is a significant departure from the ``true'' $X$-module.
This is the content of the next example.

\begin{exmp}\label{ex:comparison-percent-interleavings}
     The bound in \cref{ex:disc-direction-bound} is significantly stronger than the bound in \cref{ex: disk_eu}. This can be seen by noting that the proportion $\frac{2\pi}{6}$ in the interval $[-\pi, 0]$ is much smaller than the proportion $2\sqrt{2 - \sqrt{3}}$ in the interval $[-1, 1]$.
\end{exmp}

\subsection{Comparing Euclidean and Angular Coordinates}

A direct comparison of the error bounds for the Euclidean ($\bI$-indexed) and angular ($\Theta$-indexed) versions of the PHT is ill-defined, as their interleaving distances are computed over different indexing categories. To facilitate a meaningful comparison, we can re-parameterize the angular PHT to the Euclidean domain. Specifically, compute the extension in angular coordinates $\overline{\mathrm{PHT}}^\Theta(M)$ and then apply a change of coordinates as described in \cref{rem: change of coords} to obtain a $\bI$-persistence module. Recall the change of coordinates mentioned earlier, $\pht^\mathbb{I}_n(M)(v)(t)=\pht_n^\Theta(M)(v)(-\arccos{t})$. We may then compute the interleaving distance for this re-parameterized $\bI$-module, and in doing so, we observe the advantage of these angular coordinates. The stability of this re-parameterized $\bI$-module is given by the following lemma, which is analogous to \cref{lem: interleaving Kan}. The proof is relegated to \cref{app: lem chnage of coords stability proof}. 


\begin{lem}[Stability of the Re-parameterized Angular PHT]\label{lem:reparam_angular_stability}
For any direction $w \in \bsd$ and homological degree $n$, there exists an $\epsilon=2\sin(d_g(v, w))$-interleaving between $U'_w = \mathrm{PHT}^\Theta_n(M)(w)\circ (-\arccos)$ and $\overline{U'}_w = \overline{\mathrm{PHT}}^\Theta_n(M)(w)\circ (-\arccos)$, where $v\in A$ is any direction where $\PHT(M)$ is sampled.
Here we are using $\bI$-interleaving coordinates, as $U'_w$ and $\overline{U'}_w$ are now treated as $\bI$-modules.
\end{lem}

\begin{rem}\label{compare the two constants}
    A comparison between $\sin(\arccos(v \cdot w))$ and $\|v - w\|$ highlights the difference between two types of interleavings (cf. case $X=\bI$ in \cref{lem: interleaving Kan} and \cref{lem:reparam_angular_stability}). The former represents the perpendicular distance between vectors $v$ and $w$ on the unit sphere, while the latter represents the Euclidean distance between them. For instance, when the angle $\theta = \pi/2$, we have $\sin(\theta) = 1$, whereas the Euclidean distance between $v$ and $w$ is $\sqrt{2}$. This demonstrates that using angular coordinates for a PHT provides a tighter stability bound. In fact, it is the tightest bound one can obtain without having any other information about the shape.
\end{rem}

\subsection{Extension from Fully Discrete Data}

In practice, the Persistent Homology Transform is often computed from data that are discrete in both directions and filtration values, i.e., not just the directions $ A\subseteq\bS^{d-1}$ but for each direction $v\in A$, a finite set of filtration values $\mathbb{T}\subset \bX$.
We can handle this by defining a new category $\mathbb{AT}$ that is a full subcategory of $\mathbb{AX}$, defined previously.  
Our problem is to extend a functor $G^X: \mathbb{AT} \to \Vect$, representing the doubly-sampled PHT, to the continuous category $\bS^{d-1}\bX$. The Kan extension formulation is identical, with $\mathbb{AT}$ as the source category:
\[\begin{tikzcd}
	\mathbb{AT} &&&& \mathrm{Vect} \\ \\
	&& {\bS^{d-1}\mathbb{X}}
	\arrow[hook, from=1-1, to=3-3]
	\arrow["G^X", from=1-1, to=1-5]
	\arrow["{\hat{G}^X}"', dashed, from=3-3, to=1-5]
\end{tikzcd}\]


An argument very similar to the proof of \cref{lem: interleaving Kan} yields the following, the strongest discrete approximation result known for the PHT.
We relegate its proof to \cref{proof:stability-discrete-Kan-ext}.

\begin{thm}[Stability of the fully discrete Kan Extension]\label{thm:stability_discrete_kan_ext}
Let $M$ be a constructible set in $\mathbb{R}^d$. Let $A$ be an $\epsilon_A$-net of $\bS^{d-1}$ and let $\mathbb{T}$ be an $\epsilon_{\bT}$-net of $X$. The PHT interleaving distance between the true PHT, $\pht^X_n(M)$, and its Kan extension from the discrete data on $\mathbb{AT}$, $\overline{\pht}^X_n(M)$, is bounded for every degree $n$ as follows:
$$d_I(\pht^X_n(M), \overline{\pht}^X_n(M)) \leq 2\epsilon_A + \epsilon_\bT.$$
\end{thm}

\bibliographystyle{naturemag}
\bibliography{ref}


\appendix

\section{Kan Extensions}\label{app:Kan}

We provide a brief review of Kan extensions. For a comprehensive treatment, see \cite{cattheory}.

\begin{defn}[Left Kan Extension]
    Let $\cC$, $\cD$, and $\cE$ be categories, and let $F: \cC \to \cE$ and $G: \cC \to \cD$ be functors. A left Kan extension of $F$ along $G$ is a pair $(\mathrm{Lan}_GF, \eta)$, where $\mathrm{Lan}_GF: \cD \to \cE$ is a functor and $\eta: F \implies \mathrm{Lan}_GF \circ G$ is a natural transformation, such that for any other pair $(K, \gamma)$ with $K: \cD \to \cE$ and $\gamma: F \implies K \circ G$, there exists a unique natural transformation $\tau: \mathrm{Lan}_GF \implies K$ making the following diagram commute:

\[\begin{tikzcd}[ampersand replacement=\&]
	\cC \&\&\&\& \cE \\
	\\
	\&\& \cD
	\arrow[""{name=0, anchor=center, inner sep=0}, "F", from=1-1, to=1-5]
	\arrow["G"', from=1-1, to=3-3]
	\arrow["{\mathrm{Lan}_GF}"', dashed, from=3-3, to=1-5]
	\arrow["\eta", shorten <=11pt, shorten >=11pt, Rightarrow, from=0, to=3-3]
\end{tikzcd}\]
    
\end{defn}

The right Kan extension is defined dually and we refer the reader to \cite{cattheory} for a complete treatment. 
We are interested in computing the pointwise Kan extensions, i.e. given any object $d\in \cD$ in the domain, how do find the corresponding object $\mathrm{Lan}_GF(d)$ in $\cE$.  
To do so, we first give the definition of a comma category and its dual. 

\begin{defn}[Comma Categories]
    For a functor $G:\cC\to\cD$ and object $d\in \cD$ the \define{under comma category}, written $G\downarrow d$, consists of those objects in the image of $G$ that sit ``under'' a target object $d$, i.e., 
    \begin{itemize}
        \item  objects are pairs $(c,f)$ where $c\in \cC$ and $f:Gc\to d$ is a morphism in $\cD$,
        \item  morphisms between $(c,f)$ and $(c',f')$ is a morphism $h:c\to c'$ in $\cC$ such that the following triangle commutes: 

    \[\begin{tikzcd}[ampersand replacement=\&]
    	\&\& d \&\& \\
        \&\& \&\& \\
        Gc \&\&\&\& {Gc'}
    	\arrow["Gh", from=3-1, to=3-5]
    	\arrow["f"', to=1-3, from=3-1]
    	\arrow["{f'}", to=1-3, from=3-5]
    \end{tikzcd}\]
        
    \end{itemize}
    Dually, the \define{over comma category}, written $d \downarrow G$, consists of those objects in the image of $G$ that sit ``over'' a target object $d$, i.e.,
     \begin{itemize}
        \item  objects are pairs $(c,f)$ where $c\in \cC$ and $f:d\to Gc$ is a morphism in $\cD$,
        \item  morphisms between $(c,f)$ and $(c',f')$ is a morphism $h:c\to c'$ in $\cC$ such that the following triangle commutes: 
            \[\begin{tikzcd}[ampersand replacement=\&]
    	Gc \&\&\&\& {Gc'} \\
    	\&\& {} \\
    	\&\& d
    	\arrow["Gh", from=1-1, to=1-5]
    	\arrow["f"', to=1-1, from=3-3]
    	\arrow["{f'}", to=1-5, from=3-3]
    \end{tikzcd}\]
    \end{itemize}
\end{defn}

Both comma categories have a projection functor to $\cC$ that sends an object $(c,f)\mapsto c$ and a morphism $h:c\to c'$ to itself.
We call the projection functor from the under comma category $U_d: G\downarrow d\to \cC$ and the projection functor from the over comma category $O_d: d\downarrow G \to \cC$.

\begin{thm}[Existence of Pointwise Kan Extensions]\label{thm:ptwise-Kan}
    For $\cC$ small and $\cE$ cocomplete, the left Kan extension $\mathrm{Lan}_G F$ exists and is defined by the map $d\mapsto \mathrm{colim} F \circ U_d $. 
\end{thm}


Finally, the following definition provides an important criterion for reducing colimit calculations to simpler (sub)categories.

\begin{defn}[Cofinality]\label{defn:cofinal}
    A functor $H: \mathcal{A} \to \mathcal{B}$ is said to be \define{cofinal} if for every object $b\in \mathcal{B}$ the over comma category $(b\downarrow H)$ is non-empty and connected.


\end{defn}

\begin{rem}
    The utility of the above definition only comes to light when one learns that the above ``topological'' stipulation on the fibers of a functor $H$ is equivalent to the statement that $H$ is cofinal if for every functor $I:\mathcal{B}\to \cC$ the induced map on colimits
    \[
    \mathrm{colim}\, I\circ H \to \mathrm{colim\,} I
    \]
is an isomorphism.
\end{rem}

\section{Relegated Proofs}

\begin{lem}\label{lem:triangle ineq for geodesic}
    For any $x\in M\in\CS(\bR^d)$ and $v,w\in \bsd$ with $\|x\|\le 1$,
    $$|\arccos{x\cdot v} - \arccos{x\cdot w}|\le \arccos{v\cdot w}.$$
\end{lem}

\begin{proof}
    If $x = 0$, then $\arccos(x \cdot v) = \arccos(x \cdot w) = \pi/2$, so the inequality holds. Let $u=x/ \|x\|$ be a unit vector. Set 
\[
\alpha=\arccos(u\cdot v),\quad \beta=\arccos(u\cdot w),\quad \gamma=\arccos(v\cdot w),
\]
so $|\alpha-\beta|\le\gamma$ by the triangle inequality for the geodesic distance on $\bS^{d-1}$.  Define $g(t)=\arccos(r\cos t)$ for $t\in [0,\pi]$, then its derivative satisfies
\[
g'(t)=\frac{r\sin t}{\sqrt{1-r^2\cos^2t}}\le1,
\]
so $|g(\alpha)-g(\beta)|\le|\alpha-\beta|\le\gamma$.  Hence, $|\arccos{x\cdot v} - \arccos{x\cdot w}|\le \arccos{v\cdot w}.$
\end{proof}

\subsection{Proof of \cref{lem: pht lipschitz}}\label{app: pht lipschitz proof}
\begin{proof}
     Define the height function in direction $v\in \bS^{d-1}$ as $h_v(x)\vcentcolon= x\cdot v. $
    By the stability theorem,
    \begin{align*}
        d_\mathrm{Int}(\mathrm{PHT}^\bI_n(M)(v),\mathrm{PHT}^\bI_n(M)(w)) &\le \|h_v-h_w\|_\infty \\
        &= \sup_{x\in M}|x\cdot v-x\cdot w| \\
        &\le \sup_{x\in M}\|x\| \|v-w\|_2 \\
        &\le 1\cdot d_2(v,w),
    \end{align*}
    since we have assumed that $\|x\| \le 1$ for all $x\in M$.
    Similarly, 
    \begin{align*}
        d_\mathrm{Int}(\mathrm{PHT}_n^\Theta(M)(v),\mathrm{PHT}_n^\Theta(M)(w)) &\le \|-\arccos \circ h_v+\arccos \circ h_w\|_\infty \\
        &= \sup_{x\in M}|-\arccos(x\cdot v)+\arccos(x\cdot w)| \\
        &\le |\arccos(v\cdot w)| \\
        &\le d_g(v,w),
    \end{align*}
    where we made use of \cref{lem:triangle ineq for geodesic}. 
\end{proof}

\begin{lem}[Reparameterization of Interleavings]\label{lem:reparam_interleaving}
    Let $U, V: \Theta \to \Vect$ be two persistence modules over the indexing category $\Theta = [-\pi, 0]$ that are $\epsilon$-interleaved. Let their re-parameterizations to the indexing category $\bI = [-1,1]$ be defined by $U'(t) = U(-\arccos t)$ and $V'(t) = V(-\arccos t)$. Then the modules $U'$ and $V'$ are $\delta$-interleaved, where $\delta = 2\sin(\epsilon/2)$.
\end{lem}

\begin{proof}

An $\epsilon$-interleaving between $U$ and $V$ consists of a pair of natural transformations, $\Phi: U \Rightarrow VT_\epsilon$ and $\Psi: V \Rightarrow UT_\epsilon$, satisfying the necessary commutative diagrams. We seek the minimal $\delta \ge 0$ for which a $\delta$-interleaving $(\Phi', \Psi')$ exists for the re-parameterized modules $U'$ and $V'$.
Under the re-parameterization $t=\cos\theta$, the morphism $\Phi$ induces maps $\phi'_t: U'(t) \to V'(t')$, where $t' \vcentcolon= \cos(T_\epsilon(-\arccos t))$. The natural transformation $\Phi': U' \Rightarrow V'T_\delta$ is constructed by composing these maps with the internal morphisms of $V'$. This construction is well-defined if and only if $t' \le T_\delta(t)$ for all $t \in \bI$, which requires that $\delta \ge \sup_{t \in \bI} (t'-t)$. The commutativity conditions for the $(\Phi', \Psi')$ interleaving are satisfied because these transformations are constructed functorially from $(\Phi, \Psi)$, which already form an interleaving. The problem thus reduces to the optimization:
$$ \delta = \sup_{t \in [-1,1]} \left[ \cos(T_\epsilon(-\arccos{t})) - t \right]. $$
Letting $\theta = -\arccos t$, we seek to maximize $f(\theta) = \cos(T_\epsilon(\theta)) - \cos\theta$ for $\theta \in [-\pi,0]$. To evaluate this supremum, we must consider the definition of the shift functor on the bounded domain $\Theta$, namely $T_\epsilon(\theta) = \min(\theta+\epsilon, 0)$. We have two cases.

First, consider the domain $\theta \in [-\pi, -\epsilon]$, where $T_\epsilon(\theta) = \theta+\epsilon$. The function to maximize is $f(\theta) = \cos(\theta+\epsilon) - \cos\theta$. This quantity represents the horizontal distance between the points on the unit circle at angles $\theta$ and $\theta+\epsilon$. This distance is maximized when the chord connecting these points is horizontal, which requires the angles to be symmetric about $-\pi/2$. This configuration yields a maximal value of $2\sin(\epsilon/2)$ by the calculation below: 
\begin{align*}
    \delta = \cos\left(-\frac{\pi}{2}+\frac{\epsilon}{2}\right) - \cos\left(-\frac{\pi}{2}-\frac{\epsilon}{2}\right) 
    = \sin(\epsilon/2) - (-\sin(\epsilon/2)) = 2\sin(\epsilon/2).
\end{align*}

Second, consider the domain $\theta \in (-\epsilon, 0]$, where $T_\epsilon(\theta) = 0$. The function becomes $f(\theta) = \cos(0) - \cos\theta = 1-\cos\theta$. On $[-\pi,0]$, this function is decreasing, so its supremum on this sub-interval occurs as $\theta$ approaches the infimum of the domain, $-\epsilon$. This yields a value of $1-\cos(-\epsilon) = 1-\cos\epsilon$.

To find the global maximum, we compare the results from the two cases. Using the identity $1-\cos\epsilon = 2\sin^2(\epsilon/2)$, we see that $2\sin(\epsilon/2) \ge 2\sin^2(\epsilon/2)$ since $\sin(\epsilon/2) \in [0,1]$. The global supremum for $f(\theta)$ is therefore $2\sin(\epsilon/2)$. This choice of $\delta$ satisfies the required inequality, thereby establishing the interleaving distance.
\end{proof}

\subsection{Proof of \cref{lem:reparam_angular_stability}}\label{app: lem chnage of coords stability proof}

\begin{proof}
Let $U_w = \PHT^\Theta_n(M)(w)$ and $\overline{U}_w = \overline{\pht}^\Theta_n(M)(w)$ denote the original persistence module and the extended persistence module over the indexing category $\Theta = [-\pi, 0]$. By Lemma~\ref{lem: interleaving Kan}, for any $v \in A$, the modules $U_w$ and $\overline{U}_w$ are $\epsilon_g$-interleaved, where $\epsilon_g = 2d_g(v, w)$.

The modules in the statement of the lemma, $U'_w$ and $\overline{U'}_w$, are the re-parameterizations of $U_w$ and $\overline{U}_w$ to the indexing category $\bI = [-1,1]$. We have established in Lemma~\ref{lem:reparam_interleaving} that if two $\Theta$-persistence modules are $\epsilon_g$-interleaved, their re-parameterized counterparts in the $\bI$-domain are $\delta$-interleaved, where $\delta = 2\sin(\epsilon_g/2)$.
Substituting $\epsilon_g = 2d_g(v,w)$ into the reparameterization formula yields the interleaving distance:
$$ \epsilon = 2\sin\left(\frac{2d_g(v, w)}{2}\right) = 2\sin(d_g(v, w)). $$
This completes the proof.
\end{proof}

\subsection{Proof of \cref{thm:stability_discrete_kan_ext}}\label{proof:stability-discrete-Kan-ext}

\begin{proof}
    In line with the simplifying assumptions from \cref{rem:simplify-interleaving-proof}, we denote the true persistence module in direction $w$ as $H_w := \pht^X_n(M)(w)$ and the fully discrete Kan extension as $C_w := \overline{\pht}^X_n(M)(w)$. We note that $C_w(\beta)$ is given by the colimit of $\pht_n^X(M)$ over the discrete past light cone $L^{-}_{\mathbb{A}\mathbb{T}}(w,\beta)$, which consists of pairs $(v,t)\in A\times \mathbb{T}$ such that $t+d(v,w) \le \beta$. Our goal is to show a $(2\epsilon_A + \epsilon_\bT)$-interleaving between $H_w$ and $C_w$. Let $\epsilon = 2\epsilon_A + \epsilon_\bT$. We must construct two natural transformations
\[
\Phi: C_w \Rightarrow H_wT_{\epsilon} \qquad \text{and} \qquad \Psi: H_w \Rightarrow C_wT_{\epsilon}
\]
that satisfy the interleaving conditions.

The construction of $\Phi: C_w \Rightarrow H_w T_{\epsilon}$ is straightforward. By definition, $C_w(\beta)$ is the colimit over a diagram whose objects are pairs $(v,t) \in A \times \bT$ in the discrete past light cone of $(w,\beta)$. For any such object, the inequality $t+d(v,w) \le \beta$ holds. Since $\epsilon \ge 0$, this implies $t+d(v,w) \le T_\epsilon(\beta)$. This latter inequality defines a morphism in the full space-time category, which induces a map from the homology group $H_v(t)$ to $H_w(T_\epsilon(\beta))$. Consequently, $H_w(T_\epsilon(\beta))$ forms a co-cone over the diagram for $C_w(\beta)$. The universal property of the colimit then guarantees the existence of a unique map $\Phi_\beta: C_w(\beta) \to H_w(T_\epsilon(\beta))$, which defines the natural transformation $\Phi$.

As in the proof of \cref{lem: interleaving Kan}, the construction of $\Psi: H_w \Rightarrow C_w T_{\epsilon}$ is more involved. We define the map $\Psi_\beta: H_w(\beta) \to C_w(T_\epsilon(\beta))$ for each $\beta$ as a composition of maps arising from the discretizations. First, since $A$ is an $\epsilon_A$-net, we choose a $v \in A$ with $d(v,w) \le \epsilon_A$. The geometric inclusion $M_{w,\beta} \subseteq M_{v,\beta+\epsilon_A}$ induces a map on homology, $H_w(\beta) \to H_v(\beta+\epsilon_A)$. Second, since the filtration value $\beta+\epsilon_A$ may not be in the discrete set $\bT$, we use the fact that $\bT$ is an $\epsilon_\bT$-net to choose a $t' \in \bT$ where $\beta+\epsilon_A \le t' \le \beta+\epsilon_A+\epsilon_\bT$. This gives a second map, $H_v(\beta+\epsilon_A) \to H_v(t')$. The term $\epsilon_\bT$-net in this context means that for any value $y$ in the parameter space, there exists a point $t' \in \bT$ in the interval $[y, y + \epsilon_\bT]$.

The resulting pair $(v,t')$ is now an element of the discrete sample space $A \times \bT$.  To define a map from $H_v(t')$ into the colimit $C_w(T_\epsilon(\beta))$, we must verify that $(v,t')$ is an object in the diagram defining this colimit. The following calculation reveals this.
$$ t' + d(v,w) \le (\beta+\epsilon_A+\epsilon_\bT) + \epsilon_A = \beta + 2\epsilon_A + \epsilon_\bT = T_\epsilon(\beta). $$
So, a canonical map from $H_v(t')$ into the colimit exists and the composition of these three maps defines $\Psi_\beta$. The functoriality of this construction yields the natural transformation $\Psi$. 

Finally, verifying that the compositions $(\Psi T_\epsilon)\Phi$ and $(\Phi T_\epsilon)\Psi$ commute with the canonical shift morphisms follows from the functoriality of homology and the universal property of the colimit. Since the bound holds for any direction $w \in \bS^{d-1}$, we conclude that $d_I(\pht^X_n(M), \overline{\pht}^X_n(M)) \le 2\epsilon_A + \epsilon_\bT$.
\end{proof}

\end{document}